\documentclass[a4paper]{amsart}
\usepackage{amsthm,amsmath,amssymb,amscd,multirow}
\usepackage[latin1]{inputenc}
 \newtheorem{thm}{Theorem}
 \newtheorem{prop}[thm]{Proposition}
 \newtheorem{lemma}[thm]{Lemma}
 \newtheorem{kor}[thm]{Corollary}
 
 \theoremstyle{definition}
 \newtheorem{definition}[thm]{Definition}
 \newtheorem{ex}[thm]{Example}

 \theoremstyle{remark}
 \newtheorem{remark}[thm]{Remark}
 
 \numberwithin{equation}{section}
 \numberwithin{thm}{section}

 \def\vol{{\rm vol}}
 
 \def\Stab{{\rm Stab}}

 \def\defe{{\rm def}}
 
 \def\id{{\rm id}}
 \def\rk{{\rm rk}}

 \def\det{{\rm det}}

 \def\inv{{\rm inv}}
 
 \def\dom{{\rm dom}}

\begin{document}
\begin{title}
{The dimension of some affine Deligne-Lusztig varieties} 
\end{title}
\author{Eva Viehmann}
\address{Mathematisches Institut der Universit\"{a}t Bonn\\ Beringstrasse 1\\53115 Bonn\\Germany}
\date{}
\begin{abstract}{We prove Rapoport's dimension conjecture for affine Deligne-Lusztig varieties for $GL_h$ and superbasic $b$. From this case the general dimension formula for affine Deligne-Lusztig varieties for special maximal compact subgroups of split groups follows, as was shown in a recent paper by G\"{o}rtz, Haines, Kottwitz, and Reuman.}
\end{abstract}
\maketitle
\section{Introduction}
Let $k$ be a finite field with $q=p^r$ elements and let $\overline{k}$ be an algebraic closure. Let $F=k((t))$ and let $L=\overline{k}((t))$. Let $\mathcal{O}_F$ and $\mathcal{O}_L$ be the valuation rings. We denote by $\sigma:x\mapsto x^q$ the Frobenius of $\overline{k}$ over $k$ and also of $L$ over $F$.

Let $G$ be a split connected reductive group over $k$. Let $A$ be a split maximal torus of $G$ and $W$ the Weyl group of $A$ in $G$. For $\mu\in X_*(A)$ let $t^{\mu}$ be the image of $t\in \mathbb{G}_m(F)$ under the homomorphism $\mu:\mathbb{G}_m\rightarrow A$. Let $B$ be a Borel subgroup of $G$ containing $A$. We write $\mu_{\dom}$ for the dominant element in the orbit of $\mu\in X_*(A)$ under the Weyl group of $A$ in $G$.

We recall the definitions of affine Deligne-Lusztig varieties from \cite{Rapoport1}, \cite{GHKR}. Let $K=G(\mathcal{O}_L)$ and let $X=G(L)/K$ be the affine Grassmannian. The Cartan decomposition shows that $G(L)$ is the disjoint union of the sets $Kt^{\mu}K$ where $\mu\in X_*(A)$ is a dominant coweight. For an element $b\in G(L)$ and dominant $\mu\in X_*(A)$, the affine Deligne-Lusztig variety $X_{\mu}(b)$ is the locally closed reduced $\overline{k}$-subscheme of $X$ defined by
\begin{equation*}
X_{\mu}(b)(\overline{k})=\{g\in G(L)/K\mid g^{-1}b\sigma(g)\in Kt^{\mu}K\}.
\end{equation*}  
Left multiplication by $g\in G(L)$ induces an isomorphism between $X_{\mu}(b)$ and $X_{\mu}(gb\sigma(g)^{-1})$. Thus the isomorphism class of the affine Deligne-Lusztig variety only depends on the $\sigma$-conjugacy class of $b$.

There is an algebraic group over $F$ associated to $G$ and $b$ whose $R$-valued points (for any $F$-algebra $R$) are given by $$J(R)=\{g\in G(R\otimes_F L)\mid g^{-1}b\sigma(g)=b\}.$$ There is a canonical $J(F)$-action on $X_{\mu}(b)$.

Let $\rho$ be the half-sum of the positive roots of $G$. By $\rk_F$ we denote the dimension of a maximal $F$-split subtorus. Let $\defe_G(b)= \rk_F G - \rk_F J$. Let $\nu\in X_*(A)_{\mathbb{Q}}$ be the Newton point of $b$, compare \cite{Kottwitz}. For nonempty affine Deligne-Lusztig varieties the dimension is given by the following formula. Note that there is a simple criterion by Kottwitz and Rapoport (see \cite{KottRapo}) to decide whether an affine Deligne-Lusztig variety is nonempty.

\begin{thm}\label{thm1} Assume that $X_{\mu}(b)$ is nonempty. Then
$$\dim(X_{\mu}(b))=\langle \rho,\mu-\nu\rangle -\frac{1}{2}\defe_G(b).$$
\end{thm}

Rapoport conjectured this in \cite{Rapoport}, Conjecture 5.10 in a different form. For the reformulation compare \cite{Kottwitz2}. In \cite{Reuman2}, Reuman verifies the formula for some small groups and $b=1$. For $G=GL_n$, minuscule $\mu$ and over $\mathbb{Q}_p$ rather than over a function field, the Deligne-Lusztig varieties have an interpretation as reduced subschemes of moduli spaces of $p$-divisible groups. In this case, the corresponding dimension formula is shown by de Jong and Oort (see \cite{deJongOort}) if $b\sigma$ is superbasic and in \cite{modpdiv} for general $b\sigma$. In \cite{GHKR} 2.15, G\"{o}rtz, Haines, Kottwitz, and Reuman prove Theorem \ref{thm1} for all $b\in A(L)$. They also show in 5.8 that if there is a Levi subgroup $M$ of $G$ such that $b\in M(L)$ is basic in $M$ and if the formula is true for $M,b$ and $\mu_M$ in a certain subset of the set of all $M$-dominant coweights, then it is also true for $(G,b,\mu)$. Thus it is enough to consider superbasic elements $b$, that is elements for which no $\sigma$-conjugate is contained in a proper Levi subgroup of $G$. They show in 5.9 that it is enough to consider the case that $G=GL_h$ for some $h$ and that $b$ is basic with $m=v_t(\det (b))$ prime to $h$. In this paper we prove Theorem \ref{thm1} for this remaining case.

The strategy of the proof is as follows: We associate to the elements of $X_{\mu}(b)$ discrete invariants which we call extended semi-modules. This induces a decomposition of each connected component of $X_{\mu}(b)$ into finitely many locally closed subschemes. Their dimensions can be written as a combinatorial expression which only depends on the extended semi-module. By estimating these expressions we obtain the desired dimension formula.

For minuscule $\mu$, and over $\mathbb{Q}_p$, the group $J(\mathbb{Q}_p)$ acts transitively on the set of irreducible components of $X_{\mu}(b)$. As an application of the proof we show that for non-minuscule $\mu$, the action of $J(F)$ on this set may have more than one orbit.\\

\noindent{\it Acknowledgements.} I am grateful to M.~Rapoport for his encouragement and helpful comments. I thank R.~Kottwitz and Ng\^{o} B.-C. for their interest in my work. This work was completed during a stay at the Universit\'{e} Paris-Sud at Orsay which was supported by a fellowship within the Postdoc-Program of the German Academic Exchange Service (DAAD). I want to thank the Universit\'{e} Paris-Sud for its hospitality.

\section{Notation and conventions}\label{section2}

From now on we use the following notation: Let $G=GL_h$ and let $A$ be the diagonal torus. Let $B$ be the Borel subgroup of lower triangular matrices. For $\mu,\mu'\in X_*(A)_{\mathbb{Q}}$ we say that $\mu\preceq \mu'$ if $\mu'-\mu$ is a non-negative linear combination of positive coroots. As we may identify $X_*(A)_{\mathbb{Q}}$ with $\mathbb{Q}^h$, this induces a partial ordering on the latter set. An element $\mu=(\mu_1,\dotsc,\mu_h)\in X_*(A)\cong \mathbb{Z}^h$ is dominant if $\mu_1\leq\dotsm\leq\mu_h$.

Let $N=L^h$ and let $M_0\subset N$ be the lattice generated by the standard basis $e_0,\dotsc,e_{h-1}$. Then $K=GL_h(\mathcal{O}_L)=\Stab (M_0)$ and $g\mapsto gM_0$ defines a bijection
\begin{equation}\label{gldefxmub2}
X_{\mu}(b)(\overline{k})\cong\{M\subset N \text{ lattice }\mid \inv(M,b\sigma(M))=t^{\mu}\}.
\end{equation}
We define the volume of $M=gM_0\in X_{\mu}(b)$ to be $v_t(\det (g))$. 
 
We assume $b$ to be superbasic. The Newton point $\nu\in X_*(A)_{\mathbb{Q}}\cong\mathbb{Q}^h$ of $b$ is then of the form $\nu=(\frac{m}{h},\dotsc,\frac{m}{h})\in\mathbb{Q}^h$ with $(m,h)=1$. For $i\in\mathbb{Z}$ define $e_i$ by $e_{i+h}=te_i$. We choose $b$ to be the representative of its $\sigma$-conjugacy class that maps $e_i$ to $e_{i+m}$ for all $i$. For superbasic $b$, the condition that the affine Deligne-Lusztig variety is nonempty, namely $\nu\preceq\mu$, is equivalent to $\sum\mu_i=m$. From now on we assume this. 

For each central $\alpha\in X_*(A)$ there is the trivial isomorphism 
\begin{equation*}
X_{\mu}(b)\rightarrow X_{\mu+\alpha}(t^{\alpha}b).
\end{equation*}
We may therefore assume that all $\mu_i$ are nonnegative. For the lattices in (\ref{gldefxmub2}), this implies that $b\sigma(M)\subseteq M$.

In the following we will abbreviate the right hand side of the dimension formula for $X_{\mu}(b)$ by $d(b,\mu)$.

The set of connected components of $X$ is isomorphic to $\mathbb{Z}$, an isomorphism is given by mapping $g\in GL_h(L)$ to $v_t(\det(g))$. Let $X_{\mu}(b)^i$ be the intersection of the affine Deligne-Lusztig variety with the $i$-th connected component of $X$. Let $\pi\in GL_h(L)$ with $\pi(e_i)=e_{i+1}$ for all $i\in\mathbb{Z}$. Then $\pi$ commutes with $b\sigma$, and defines isomorphisms $X_{\mu}(b)^i\rightarrow X_{\mu}(b)^{i+1}$ for all $i$. Thus it is enough to determine the dimension of $X_{\mu}(b)^0$.

For superbasic $b$, an element of $J(F)$ is determined by its value at $e_0$. More precisely, $J(F)$ is the multiplicative subgroup of a central simple algebra over $F$. Hence $\defe_G(b)=h-1$. If $v_t(\det (g))=i$ for some $g\in J(F)$, then $g$ induces isomorphisms between $X_{\mu}(b)^{j}$ and $X_{\mu}(b)^{j+i}$ for all $j$. On $X_{\mu}(b)^0$, we have an action of $\{g\in J(F)\mid v_t(\det(g))=0\}=J(F)\cap \Stab(M_0)$.

\begin{remark}\label{remgp}
To a vector $\psi=(\psi_i)\in\mathbb{Q}^h$ we associate the polygon in $\mathbb{R}^2$ that is the graph of the piecewise linear continuous function $f:[0,h]\rightarrow\mathbb{R}$ with $f(0)=0$ and slope $\psi_i$ on $[i-1,i]$. One can easily see that $d(b,\mu)$ is equal to the number of lattice points below the polygon corresponding to $\nu$ and (strictly) above the polygon corresponding to $\mu$. 
\end{remark}
 
\section{Extended semi-modules}
In this section we describe the combinatorial invariants which are used to decompose $X_{\mu}(b)^0$.

\begin{definition}\label{defsemimod}
\begin{enumerate}
\item Let $m$ and $h$ be coprime positive integers. A \emph{semi-module} for $m$, $h$ is a subset $A\subset \mathbb{Z}$ that is bounded below and satisfies $m+A\subset A$ and $h+A\subset A$. Let $B=A\setminus (h+A)$. The semi-module is called \emph{normalized} if $\sum_{a\in B}a=\frac{h(h-1)}{2}$.
\item Let $\nu=(\frac{m}{h},\dotsc,\frac{m}{h})\in\mathbb{Q}^h$. Let $\mu'=(\mu'_1,\cdots,\mu'_h)\in \mathbb{N}^h$ not necessarily dominant with $\nu\preceq \mu'$. A semi-module $A$ for $m$, $h$ is \emph{of type $\mu'$} if the following condition holds: Let $b_0=\min\{b\in B\}$ and let inductively $b_{i}=b_{i-1}+m-\mu'_ih\in\mathbb{Z}$ for $i=1,\dotsc,h$. Then $b_0=b_h$ and $\{b_i\mid i=0,\dotsc,h-1 \}=B$. 
\end{enumerate}
\end{definition}
\begin{remark}
Semi-modules are also used by de Jong and Oort in \cite{deJongOort} to define a stratification of a moduli space of $p$-divisible groups whose rational Dieudonn\'{e} modules are simple of slope $\frac{m}{h}$. In this case $\mu$ is minuscule, and they use semi-modules for $m,h-m$ to decompose the moduli space.
\end{remark}
\begin{lemma}
If $A$ is a semi-module, then its translate $-\frac{\sum_{a\in B}a }{h}+\frac{h-1}{2}+A$ is the unique normalized translate of $A$. It is called the normalization of $A$. There is a bijection between the set of normalized semi-modules for $m$, $h$ and the set of possible types $\mu'\in\mathbb{N}^h$ with $\nu\preceq\mu'$.
\end{lemma}
\begin{proof}
For the first assertion one only has to notice that the fact that the $h$ elements of $B$ are incongruent modulo $h$ implies that $\sum_{a\in B}a-\frac{h(h-1)}{2}$ is divisible by $h$. For the second assertion let $A$ be a normalized semi-module, let $b_0=\min\{a\in B\}$ and let inductively $b_{i}=b_{i-1}+m-\mu'_ih$ where $\mu'_i$ is maximal with $b_i\in A$. Then $b_h=b_0$ and $\{b_i\mid i=0,\dotsc,h-1\}=B$. From $b_0<b_{i_0}$ for $i_0=1,\dotsc,h-1$ we obtain $\sum_{i=1}^{i_0}(m-\mu'_ih)>0$ for all $i_0<h$. Similarly, $b_0=b_h$ implies $\sum_{i=1}^h\mu_i'=m$. This shows $\nu\preceq\mu'$. As $m+A\subset A$, the $\mu_i'$ are nonnegative. Given $\mu'$ as above, the corresponding normalized semi-module $A$ can be constructed as follows: Let $b_0=0$, and inductively $b_{i}=b_{i-1}+m-\mu'_ih$. Then $A$ is the normalization of $\{b_i+\alpha h\mid \alpha\in \mathbb{N}, 0\leq i<h\}$.  
\end{proof}

\begin{definition}\label{defenlsm}
Let $m$ and $h$ be as before and let $\mu=(\mu_i)\in\mathbb{N}^h$ be dominant with $\sum \mu_i=m$. An \emph{extended semi-module} $(A,\varphi)$ for $\mu$ is a normalized semi-module $A$ for $m$, $h$ together with a function $\varphi:\mathbb{Z}\rightarrow \mathbb{N}\cup \{-\infty\}$ with the following properties:
\begin{enumerate}
\item $\varphi(a)=-\infty$ if and only if $a\notin A$.
\item $\varphi(a+h)\geq \varphi(a)+1$ for all $a$.
\item $\varphi(a)\leq\max\{n\mid a+m-nh\in A\}$ for all $a\in A$. If $b\in A$ for all $b\geq a$, then the two sides are equal.
\item There is a decomposition of $A$ into a disjoint union of sequences $a_j^1,\cdots,a_j^h$ with $j\in\mathbb{N}$ and the following properties:
  \begin{enumerate}
  \item $\varphi(a_{j+1}^l)=\varphi(a_j^l)+1$
  \item If $\varphi(a_j^l+h)=\varphi(a_j^l)+1$, then $a_{j+1}^l=a_j^l+h$. Otherwise $a_{j+1}^l>a_j^l+h$.
  \item The $h$-tuple $(\varphi(a_0^l))$ is a permutation of $\mu$.
  \end{enumerate}
\end{enumerate}
An extended semi-module such that equality holds in (3) for all $a\in A$ is called \emph{cyclic}.

Let $A$ be a normalized semi-module for $m$, $h$ and let $\mu'$ be its type. Let $\mu=\mu'_{\dom}$. Let $\varphi$ be such that (1) holds and that we have equality in (3) for all $a\in A$. Then in (2) the two sides are also equal for all $a\in A$. A decomposition of $A$ as in (4) is given by putting all elements into one sequence that are congruent modulo $h$. Hence $(A,\varphi)$ is a cyclic extended semi-module for $\mu$, called the \emph{cyclic extended semi-module associated to $A$}.
\end{definition}
\begin{ex}\label{ex2}
We give an explicit example of a non-cyclic extended semimodule for $m=4$, $h=5$, and $\mu=(0,0,0,2,2)$. Let $A$ be the normalized semi-module of type $(0,0,1,2,1)$. Then $B=A\setminus (5+A)$ consists of $-2,-1,2,5$, and $6$. Let $\varphi(-1)=0$ and $\varphi(a)=\max\{n\mid a+m-nh\in A\}$ if $a\in A\setminus \{-1\}$. See also Figure \ref{pic1} that shows elements of $A$ marked by crosses and the corresponding values of $\varphi$. A decomposition of $A$ is given as follows: Three sequences are given by the elements of $A$ congruent to $-2$, $2$, and $5$ modulo $5$, respectively. The forth sequence is given by all elements congruent to $4$ modulo $5$ and greater than $-1$. The last sequence consists of the remaining elements $-1$ and $6,11,16,\dotsc$.
{\setlength{\unitlength}{12 pt}
\begin{figure}[h]
\begin{center}
\begin{picture}(25,4)(0,0.5)

\put(1.5,1){\makebox(0,0){$\varphi(a)$}}
\put(5,1){\makebox(0,0){$\dotsm$}}
\put(6.5,1){\makebox(0,0){$-\infty$}}
\put(8,1){\makebox(0,0){$0$}}
\put(9.5,1){\makebox(0,0){$0$}}
\put(11,1){\makebox(0,0){$-\infty$}}
\put(12.5,1){\makebox(0,0){$-\infty$}}
\put(14,1){\makebox(0,0){$0$}}
\put(15.5,1){\makebox(0,0){$1$}}
\put(17,1){\makebox(0,0){$2$}}
\put(18.5,1){\makebox(0,0){$2$}}
\put(20,1){\makebox(0,0){$1$}}
\put(21.5,1){\makebox(0,0){$1$}}
\put(23,1){\makebox(0,0){$2$}}
\put(24.5,1){\makebox(0,0){$\dotsm$}}

\put(1.5,4){\makebox(0,0){$a$}}
\put(6.5,4){\makebox(0,0){$-3$}}
\put(8,4){\makebox(0,0){$-2$}}
\put(9.5,4){\makebox(0,0){$-1$}}
\put(11,4){\makebox(0,0){$0$}}
\put(12.5,4){\makebox(0,0){$1$}}
\put(14,4){\makebox(0,0){$2$}}
\put(15.5,4){\makebox(0,0){$3$}}
\put(17,4){\makebox(0,0){$4$}}
\put(18.5,4){\makebox(0,0){$5$}}
\put(20,4){\makebox(0,0){$6$}}

\put(3.5,2.5){\makebox(0,0){$\dotsm$}}
\put(5,2.5){\makebox(0,0){$\cdot$}}
\put(6.5,2.5){\makebox(0,0){$\cdot$}}
\put(8,2.5){\makebox(0,0){$\times$}}
\put(9.5,2.5){\makebox(0,0){$\times$}}
\put(11,2.5){\makebox(0,0){$\cdot$}}
\put(12.5,2.5){\makebox(0,0){$\cdot$}}
\put(14,2.5){\makebox(0,0){$\times$}}
\put(15.5,2.5){\makebox(0,0){$\times$}}
\put(17,2.5){\makebox(0,0){$\times$}}
\put(18.5,2.5){\makebox(0,0){$\times$}}
\put(20,2.5){\makebox(0,0){$\times$}}
\put(21.5,2.5){\makebox(0,0){$\times$}}
\put(23,2.5){\makebox(0,0){$\times$}}
\put(24.5,2.5){\makebox(0,0){$\dotsm$}}

\end{picture}
\end{center}
\caption{A non-cyclic extended semi-module} \label{pic1}
\end{figure}
}

\end{ex}

\begin{lemma}\label{lemnoncyclic}
If $(A,\varphi)$ is an extended semi-module for $\mu$, and if $\mu^0$ is the type of $A$, then $\mu^0_{\dom}\preceq \mu$. If $\mu^0_{\dom}=\mu$, then $(A,\varphi)$ is a cyclic extended semi-module.
\end{lemma}
\begin{proof}
Let $(A,\varphi_0)$ be the cyclic extended semi-module associated to $A$. Let $$\{x_1,\dotsc,x_n\}=\{a\in A\mid \varphi(a+h)>\varphi(a)+1\}$$ with $x_i>x_{i+1}$ for all $i$. For $i\in \{1,\dotsc,n\}$ let 
\begin{equation*}
\varphi_{i}(a)=
\begin{cases}
-\infty&\text{if } a\notin A\\
\varphi(a)&\text{if } a\geq x_i\\
\varphi_i(a+h)-1&\text{else.}
\end{cases}
\end{equation*}
We show that $(A,\varphi_i)$ is an extended semi-module for some $\mu^i$ with $\mu_{\dom}^{i-1}\preceq \mu_{\dom}^{i}$ and $\mu_{\dom}^{i-1}\neq \mu_{\dom}^i$ for all $i\geq 1$. As $\varphi_n=\varphi$, it then follows that $\mu^0_{\dom}\preceq \mu^n_{\dom}=\mu$ with equality if and only if $n=0$, that is if $\varphi$ is cyclic.

The decomposition of $(A,\varphi_{i})$ is defined as follows: For $a<x_i$, the successor of $a$ is $a+h$. Otherwise it is the successor from the decomposition of $(A,\varphi)$. From the properties of the decompositions for $\varphi_0$ and $\varphi$ one deduces that the decomposition satisfies the required properties. Let $n_i\geq 0$ be maximal with $x_i-n_ih\in A$ and let $\alpha_i=\varphi(x_i+h)-1-\varphi(x_i)>0$. Thus $\varphi_i$ is obtained from $\varphi_{i-1}$ by subtracting $\alpha_i$ from the values at $x_i,x_i-h\dotsc,x_i-n_ih$. From $\mu^{i-1}$ we obtain $\mu^i$ by replacing the two entries $\varphi_{i-1}(x_i-n_ih)=\varphi_{i-1}(x_i)-n_i$ and $\varphi_{i-1}(x_i)-\alpha_i+1$ (which is the value of $\varphi$ of the successor of $x_i$ in the sequence corresponding to $\varphi_i$) by $\varphi_{i-1}(x_i)-\alpha_i-n_i$ and $\varphi_{i-1}(x_i)+1$. As $$\varphi_{i-1}(x_i)-n_i, \varphi_{i-1}(x_i)-\alpha_i+1\in (\varphi_{i-1}(x_i)-\alpha_i-n_i,\varphi_{i-1}(x_i)+1),$$ we have $\mu_{\dom}^{i-1}\preceq \mu_{\dom}^{i}$ and $\mu_{\dom}^{i-1}\neq \mu_{\dom}^i$.
\end{proof}

\begin{kor}
If $\mu$ is minuscule, then all extended semi-modules for $\mu$ are cyclic.
\end{kor}
\begin{proof}
Let $(A,\varphi)$ be such an extended semi-module. Let $\mu'$ be the type of $A$. Then $\mu'_{\dom}\preceq \mu$, thus $\mu'_{\dom}=\mu$. Hence the assertion follows from the preceding lemma. 
\end{proof}

\begin{lemma}\label{lemmufin}
There are only finitely many extended semi-modules $(A,\varphi)$ for each $\mu$.
\end{lemma}
\begin{proof}
Let $\mu'$ be the type of the semi-module $A$. As $\mu'_{\dom}\preceq \mu$, there are only finitely many possible types and corresponding normalized semi-modules. For fixed $A$, the third condition for extended semi-modules determines all but finitely many values of $\varphi$. For the remaining values we have $0\leq \varphi(a)\leq\max\{n\mid a+m-nh\in A\}$. Thus for each $A$ there are only finitely many possible functions $\varphi$ such that $(A,\varphi)$ is an extended semi-module for $\mu$.
\end{proof}

\section{The decomposition of the affine Deligne-Lusztig variety}

Let $M\in X_{\mu}(b)^0$ be a lattice in $N$. In this section we associate to $M$ an extended semi-module for $\mu$. This leads to a paving of $X_{\mu}(b)^0$ by finitely many locally closed subschemes. For minuscule $\mu$, this decomposition of the set of lattices is the same as the one constructed by de Jong and Oort in \cite{deJongOort}, compare also \cite{modpdiv}, Section 5.1.

Let $m$ and $h$ be as in Section \ref{section2}. Let $v\in N$ and recall that $te_i=e_{i+h}$. Then we can write $v=\sum_{i\in\mathbb{Z}}\alpha_ie_i$ with $\alpha_i\in\overline{k}$ and $\alpha_i=0$ for small $i$. Let
\begin{eqnarray*}
I:N\setminus \{0\}&\rightarrow&\mathbb{Z}\\
v&\mapsto&\min\{i\mid \alpha_i\neq 0\}.
\end{eqnarray*}
For a lattice $M\in X_{\mu}(b)^0$ we consider the set
\begin{equation*}
A=A(M)=\{I(v)\mid v\in M\setminus \{0\}\}.
\end{equation*}
Then $A(M)$ is bounded below and $h+A(M)\subset A(M)$. As $b\sigma(M)\subset M$, we have $m+A(M)\subset A(M)$, thus $A(M)$ is a semi-module for $m$, $h$. 
We have
\begin{equation*}
\vol(M)=\mid \mathbb{N}\setminus (A\cap\mathbb{N} )\mid - \mid A\setminus (\mathbb{N}\cap A)\mid=0.
\end{equation*}
This implies that $\sum_{a\in B}a=\sum_{i=0}^{h-1}i$, thus $A$ is normalized.

Let further
\begin{eqnarray*}
\varphi=\varphi(M):\mathbb{Z}&\rightarrow&\mathbb{N}\cup \{-\infty\}\\
a&\mapsto&\begin{cases} 
\max\{n\mid \exists v\in M \text{ with } I(v)=a, t^{-n}b\sigma(v)\in M\}&\text{if }a\in A(M)\\
-\infty&\text{else.}
\end{cases}
\end{eqnarray*}
Note that by the definition of $A(M)$, the set on the right hand side is nonempty. As $b\sigma(M)\subset M$, the values of $\varphi$ are indeed in $\mathbb{N}\cup \{-\infty\}$.
 
\begin{lemma}\label{lemphim}
Let $M\in X_{\mu}(b)^0$. Then $(A(M),\varphi(M))$ is an extended semi-module for $\mu$.
\end{lemma}
\begin{proof}
We already saw that $A(M)$ is a normalized semi-module. We have to check the conditions on $\varphi$. The first condition holds by definition. Let $v\in M$ with $I(v)=a$ be realizing the maximum for $\varphi(a)$. Then $tv\in M$ with $I(tv)=a+h$ implies that $\varphi(a+h)\geq\varphi(a)+1$, which shows (2). Let $v\in M$ with $I(v)=a$ and $t^{-\varphi(a)}b\sigma(v)\in M$. Then $I(t^{-\varphi(a)}b\sigma(v))=a+m-\varphi(a)h\in A(M)$, whence the first part of (3). Let $b\in A$ for all $b\geq a$. Let $n_0=\max \{n\mid a+m-nh\in A\}$. Let $v'\in M$ with $I(v')=a+m-n_0h$ and let $v=(b\sigma)^{-1}(t^{n_0}v')\in N$. Then $I(v)=a$, thus $v=\sum_{b\geq a}\alpha_b e_b$ for some $\alpha_b\in\overline{k}$. As $b\in A$ for all $b\geq a$, we also have $e_b\in M$ for all such $b$. Thus $v\in M$ with $t^{-n_0}b\sigma(v)=v'\in M$. Hence $\varphi(a)= n_0$. It remains to show (4). For $a\in \mathbb{Z}$ and $\varphi_0\in\mathbb{N}$ let $$\tilde{V}_{a,\varphi_0}=\{v\in M\mid v=0\text{ or }I(v)\geq a, t^{-\varphi_0}b\sigma(v)\in M\}$$ and $V_{a,\varphi_0}=\tilde{V}_{a,\varphi_0}/\tilde{V}_{a,\varphi_0+1}$. Then $V_{a,\varphi_0}$ is a $\overline{k}$-vector space of dimension $\mid \{a'\geq a\mid \varphi(a')=\varphi_0\}\mid$. We construct the sequences by inductively sorting all elements $a\in A$ with $\varphi(a)\leq \varphi_0$ for some $\varphi_0$: For $\varphi_0=\min\{\varphi(a)\mid a\in A\}$ we take each element $a$ with this value of $\varphi$ as the first element of a sequence. (At the end we will see that we did not construct more than $h$ sequences.) We now describe the induction step from $\varphi_0$ to $\varphi_0+1$: If $v_1,\dotsc,v_i$ is a basis of $V_{a,\varphi_0}$ for some $a$, then the $tv_j$ are linearly independent in $V_{a+h,\varphi_0+1}$. Thus $\dim V_{a,\varphi_0}\leq \dim V_{a+h,\varphi_0+1}$ for every $a$. Hence there are enough elements $a\in A$ with $\varphi(a)=\varphi_0+1$ to prolong all existing sequences such that conditions (a) and (b) are satisfied. We take the $a\in A$ with $\varphi(a)=\varphi_0+1$ that are not already in some sequence as first elements of new sequences. Inductively, this constructs sequences with properties (a) and (b). To show (c), let $a<b_0$. Then
\begin{eqnarray*}
\mid \{i\mid\mu_i=n \} \mid&=&\dim_{\overline{k}}V_{a,n}-\dim_{\overline{k}}V_{a-h,n-1}\\
&=&\mid \{a^l_0\mid \varphi(a^l_0)=n\}\mid. 
\end{eqnarray*} 
This also shows that we constructed exactly $h$ sequences.
\end{proof}

For each extended semi-module $(A,\varphi)$ for $\mu$ let $$\mathcal{S}_{A,\varphi}=\{M\subset N \text{ lattice}\mid A(M)=A, \varphi(M)=\varphi\}\subset X.$$

\begin{lemma}\label{lemdim}
The sets $\mathcal{S}_{A,\varphi}$ are contained in $X_{\mu}(b)^0$. They define a decomposition of $X_{\mu}(b)^0$ into finitely many disjoint locally closed subschemes. Especially, $\dim X_{\mu}(b)^0=\max\{\dim \mathcal{S}_{A,\varphi}\}$.
\end{lemma}
\begin{proof}
The last property in the definition of an extended semi-module shows that $(A,\varphi)$ determines $\mu$. Thus $\mathcal{S}_{A,\varphi} \subseteq X_{\mu}(b)^0$. Using Lemma \ref{lemmufin} and Lemma \ref{lemphim} it only remains to show that the subschemes are locally closed. The condition that $a\in A(M)$ is equivalent to $\dim (M\cap \langle e_a,e_{a+1},\dotsc\rangle)/(M\cap \langle e_{a+1},e_{a+2},\dotsc\rangle)=1$. This is clearly locally closed. If $a$ is sufficiently large, it is contained in all extended semi-modules for $\mu$ and if $a$ is sufficiently small, it is not contained in any extended semi-module for $\mu$. Thus fixing $A$ is an intersection of finitely many locally closed conditions on $X_{\mu}(b)^0$, hence locally closed. Similarly, it is enough to show that $\varphi(a)< n$ for some $a\in A$ and $n\in \mathbb{N}$ is an open condition on $\{M\in X\mid b\sigma(M)\subset M, A(M)=A\}\subset X$. But this condition is equivalent to 
$$(\langle e_i\mid i\geq a\rangle\cap M \cap t^n(b\sigma)^{-1}(M))/\langle e_i\mid i\geq a+1\rangle= (0),$$ which is an open condition.
\end{proof}

Let $(A,\varphi)$ be an extended semi-module for $\mu$. Let 
\begin{equation}
\mathcal{V}(A,\varphi)=\{(a,b)\in A\times A\mid b>a, \varphi(a)>\varphi(b)>\varphi(a-h)\}.
\end{equation} 
\begin{thm}\label{thm2} \begin{enumerate}
\item Let $A$ and $\varphi$ be as above. There exists a nonempty open subscheme $U(A,\varphi)\subseteq\mathbb{A}^{\mathcal{V}(A,\varphi)}$ and a morphism $U(A,\varphi)\rightarrow \mathcal{S}_{A,\varphi}$ that induces a bijection between the set of $\overline{k}$-valued points of $U(A,\varphi)$ and $\mathcal{S}_{A,\varphi}$. Especially, $\dim(\mathcal{S}_{A,\varphi})=\mid \mathcal{V}(A,\varphi)\mid$.
\item If $(A,\varphi)$ is a cyclic extended semi-module, then $U(A,\varphi)=\mathbb{A}^{\mathcal{V}(A,\varphi)}$.
\end{enumerate}
\end{thm}
\begin{proof}
We denote the coordinates of a point $x$ of $\mathbb{A}^{\mathcal{V}(A,\varphi)}$ by $x_{a,b}$ with $(a,b)\in \mathcal{V}(A,\varphi)$. To define a morphism $\mathbb{A}^{\mathcal{V}(A,\varphi)}\rightarrow X$, we describe the image $M(x)$ of a point $x\in\mathbb{A}^{\mathcal{V}(A,\varphi)}(R)$ where $R$ is a $\overline{k}$-algebra. For each $a\in A$ we define an element $v(a)\in N_R=N\otimes_{\overline{k}}R$ of the form $v(a)=\sum_{b\geq a}\alpha_b e_b$ with $\alpha_a=1$. The $R[[t]]$-module $M(x)\subset N_R$ will then be generated by the $v(a)$. We want the $v(a)$ to satisfy the following relations: For $a\in h+A$ we want
\begin{equation}\label{glrec2}
v(a)=tv(a-h)+\sum_{(a,b)\in\mathcal{V}(A,\varphi)}x_{a,b}v(b).
\end{equation}
Let $y=\max\{b\in B\}$. If $a=y$ we want 
\begin{equation}\label{glrec30}
v(a)=e_a+\sum_{(a,b)\in \mathcal{V}(A,\varphi)}x_{a,b}v(b).
\end{equation}
For all other elements $a\in B$, we want the following equation to hold: Let $a'\in A$ be minimal with $a'+m-\varphi(a')h=a$. Then $v'=t^{-\varphi(a')}b\sigma(v(a'))\in N_R$ with $I(v')=a$. Let
\begin{equation}\label{glrec3}
v(a)=v'+\sum_{(a,b)\in \mathcal{V}(A,\varphi)}x_{a,b}v(b).
\end{equation}

{\bf Claim 1.} For every $x\in \mathbb{A}^{\mathcal{V}(A,\varphi)}(R)$ there are uniquely determined $v(a)\in N_R$ for all $a\in A$ satisfying (\ref{glrec2}) to (\ref{glrec3}).

We set 
\begin{equation*}
v(a)=\sum_{j\in\mathbb{N}}\alpha_{a,j}e_{a+j}
\end{equation*}
with $\alpha_{a,j}\in R$ and $\alpha_{a,0}=1$ for all $a$. We solve the equations by induction on $j$. Assume that the $\alpha_{a,j}$ are determined for $j\leq j_0$ and such that the equations for $v(a)$ hold up to summands of the form $\beta_je_j$ with $j>a+j_0$. To determine the $\alpha_{a,j_0+1}$, we write $a\equiv y+im \pmod{h}$ and proceed by induction on $i\in \{0,\dotsc,h-1\}$. For $i=0$ and $a=y$, the coefficient $\alpha_{a,j_0+1}$ is the uniquely determined element such that (\ref{glrec30}) holds up to summands of the form $\beta_{j}e_j$ with $j>j_0+1$. Note that by induction on $j$ and as $b>a$, the coefficient of $e_{y+j_0+1}$ on the right hand side of the equation is determined. For $a=y+nh$ with $n>0$, the coefficients are similarly defined by (\ref{glrec2}). For $i>0$ and $a\in A$ minimal in this congruence class, the coefficient is determined by (\ref{glrec3}). Here, the coefficient of $e_{a+j_0+1}$ on the right hand side of each equation is determined by induction on $i$ and $j$. For larger $a$ in this congruence class we use again (\ref{glrec2}). By passing to the limit on $j$, we obtain the uniquely defined $v(a)\in N_R$ solving the equations.

{\bf Claim 2.} Let $M(x)=\langle v(a)\mid a\in A\rangle_{R[[t]]}$. Then at each specialization of $x$ to a $\overline{k}$-valued point $y$ we have $A= A(M(y))$ and $\varphi(M(y))(a)\geq\varphi(a)$ for all $a$. 

From the definition of $M$ we immediately obtain $A\subseteq A(M(y))$. To show equality consider an element $v=\sum_a \alpha_av(a)\in M(y)=M$. Write $v=\sum_{i\in\mathbb{Z}} b_ie_i$ with $b_i\in \overline{k}$. Let $i_0=\min\{I(\alpha_av(a))\}$. If $b_{i_0}\neq 0$, then $I(v)=i_0\in A$. Otherwise we consider $\sum_{\{a\mid I(\alpha_av(a))=i_0\}}\alpha_av(a)$. Note that $I(v(a))\equiv i_0\pmod{h}$ for all $a$ occuring in the sum. Then (\ref{glrec2}) shows that this sum can be written as a sum of $v(b)$ with $b>i_0$. Thus we may replace $i_0$ by a larger number. As $i\in A$ for all sufficiently large $i$, this shows that $I(v)\in A$, so $A(M)=A$.

Let $x\in \mathbb{A}^{\mathcal{V}(A,\varphi)}(\overline{k})$ and let $M=M(x)$. We show that $t^{-\varphi(a)}b\sigma(v(a))\in M$ for all $a$. This means that $\varphi(M)(a)\geq\varphi(a)$ for all $a$. Consider the elements $a'\in A$ that are minimal with $a'+m-\varphi(a')h=a$ for some $a\in B\setminus \{y\}$. For these elements, the assertion follows from (\ref{glrec3}). If $a$ is minimal with $a+m-\varphi(a)h=y$, then $I(t^{-\varphi(a)}b\sigma(v(a)))=y$. As all $e_i$ with $i\geq y$ are in $M$, this element is also contained in $M$. If $\varphi(a)=\varphi(a-h)+1$ then $v(a)=tv(a-h)$ and the assertion holds for $a-h$ if and only if it holds for $h$. From this, we obtain the claim for all $a\in A$ with $\varphi(a)=\max\{n\mid a+m-nh\in A\}$. Especially, it follows for all sufficiently large elements of $A$. It remains to prove the claim for the finitely many elements $a\in A$ with $\max\{n\mid a+m-nh\in A\}>\varphi(a)$. We use decreasing induction on $a$: Let $a$ be in this set, and assume that we know the assertion for all $a'>a$. From (\ref{glrec2}) we obtain that 
\begin{eqnarray*}
t^{-\varphi(a)}b\sigma(v(a))&=& t^{-\varphi(a)-1}b\sigma(tv(a))\\
&=& t^{-\varphi(a)-1}b\sigma(v(a+h)-\sum_{b>a+h,\varphi(a+h)>\varphi(b)\geq \varphi(a)+1}x_{a+h,b}v(b)).
\end{eqnarray*}
By induction, the right hand side is in $M$ and Claim 2 is shown. 

As all $\mu_i$ are nonnegative, we constructed a morphism from $\mathbb{A}^{\mathcal{V}(A,\varphi)}$ to the subscheme $X_A$ of $X$ defined by $X_A(\overline{k})=\{M\mid A(M)=A,b\sigma(M)\subseteq M\}$. 

{\bf Claim 3.} There is a nonempty open subscheme $U(A,\varphi)$ of $\mathbb{A}^{\mathcal{V}(A,\varphi)}$ that is mapped to $\mathcal{S}_{A,\varphi}$. If $(A,\varphi)$ is cyclic, then $U(A,\varphi)=\mathbb{A}^{\mathcal{V}(A,\varphi)}$.

In general we do not have $\varphi(M)(a)=\varphi(a)$ for all $a$. The proof of Lemma \ref{lemdim} shows that $\varphi(M)(a)\leq \varphi(a)$ is an open condition on $X_A$, and thus on $\mathbb{A}^{\mathcal{V}(A,\varphi)}$. Let $U(A,\varphi)$ be the corresponding open subscheme, which is then mapped to $\mathcal{S}_{A,\varphi}$. We have to show that it is nonempty, thus to construct a point in $\mathbb{A}^{\mathcal{V}(A,\varphi)}$ where the corresponding function $\varphi(M)$ is equal to $\varphi$. If $\varphi(a)=\max\{n\mid a+m-nh\in A\}$, then $\varphi(M)(a)=\varphi(a)$. Especially, the two functions are equal for all $a$ if $(A,\varphi)$ is cyclic. In this case $U(A,\varphi)=\mathbb{A}^{\mathcal{V}(A,\varphi)}$. If $\varphi(a)+1=\varphi(a+h)$ and if $\varphi(M)(a+h)=\varphi(a+h)$, then $\varphi(M)(a+h)-1\geq\varphi(M)(a)\geq \varphi(a)$ implies that $\varphi(M)(a)=\varphi(a)$. Thus it is enough to find a point where $\varphi(M)(a)=\varphi(a)$ for all $a\in A$ with $\varphi(a+h)>\varphi(a)+1$. For each such $a$ let $b_a$ be the successor in a decomposition of $(A,\varphi)$ into sequences. Then $(a+h,b_a)\in\mathcal{V}(A,\varphi)$. Let $x_{a+h,b_a}=1$ for these pairs and choose all other coefficients to be $0$. Then for this point and $a$ as before we have that $\varphi(M)(a)=\varphi(b_a)-1=\varphi(a)$. Thus $U(A,\varphi)$ is nonempty.

{\bf Claim 4.} The map $U(A,\varphi)\rightarrow \mathcal{S}_{A,\varphi}$ defines a bijection on $\overline{k}$-valued points.

More precisely, we have to show that for each $M\in \mathcal{S}_{A,\varphi}$ there is exactly one $x\in U(A,\varphi)(\overline{k})$ such that $M$ contains a set of elements $v(a)$ for $a\in A$ with $I(v(a))=a$ and satisfying (\ref{glrec2}) to (\ref{glrec3}) for this $x$. The argument is similar as the construction of $v(a)$ for given $x$: By induction on $j$ we will show the following assertion: There exist $x^j=(x^j_{a,b})\in U(A,\varphi)(\overline{k})$ and $v_j(a)\in M$ for all $a$ with $t^{-\varphi(a)}b\sigma(v_j(a))\in M$ and which satisfy equations (\ref{glrec2}) to (\ref{glrec3}) for $x^j$ up to summands of the form $\beta_ne_n$ with $n>a+j$. Furthermore the $x^j_{a,b}$ with $b-a\leq j$ and the coefficients of $e_n$ in $v_j(a)$ for $n\leq a+j$ will be chosen independently of $j$ and only depending on $M$.

For $j=0$ choose any $x^0\in U(A,\varphi)(\overline{k})$ and $v_0(a)\in M$ with $I(v_0(a))=a$, first coefficient $1$ and $t^{-\varphi(a)}b\sigma(v_0(a))\in M$. The existence of these $v_0(a)$ follows from $M\in X_{\mu}(b)$. Assume that the assertion is true for some $j_0$. For $n\leq j_0$ let $x^{j_0+1}_{a,a+n}=x^{j_0}_{a,a+n}$. We proceed again by induction on $i$ to define the coefficients for $a\equiv y+im \pmod{h}$. Let $a=y$.  Choose the coefficients $x^{j_0+1}_{y,y+n}$ with $n>j_0$ such that
$$v_{j_0+1}(y)=e_{y}+\sum_{(y,y+n)\in\mathcal{V}(A,\varphi)}x^{j_0+1}_{y,y+n}v_{j_0}(y+n)$$ satisfies $t^{-\varphi(y)}b\sigma(v_{j_0+1}(y))\in M$. The definition of $\varphi=\varphi(M)$ shows that such coefficients exist and from $\varphi(y+n)<\varphi(y)$ it follows that they are unique. For the other elements $v(a)$ we proceed similarly: For those with $a-h\notin A$ we use equation (\ref{glrec3}), on the right hand side with the values from the induction hypothesis, to define the new $v_{j_0+1}(a)$. For $a\in h+A$ we use (\ref{glrec2}). As we know that $t^{-\varphi(a-h)-1}b\sigma(tv_{j_0}(a-h))\in M$, it is sufficient to consider the $b>a$ with $\varphi(a-h)<\varphi(b)<\varphi(a)$. At each step the coefficient of $e_{a+j_0+1}$ of the right hand side is already defined by the induction hypothesis. It only depends on the $x_{a,a+n}^{j_0}$ and the coefficients of $e_{b+n}$ of $v_{j_0}(b)$ with $n\leq j_0$, hence only on $M$. The coefficients of $x^{j_0+1}$ are given by requiring that $t^{-\varphi(a)}b\sigma(v_{j_0+1}(a))\in M$.
\end{proof}

\section{Combinatorics}

In this section we estimate $\mid\mathcal{V}(A,\varphi)\mid$ to determine the dimension of the affine Deligne-Lusztig variety $X_{\mu}(b)$.
\begin{remark}\label{remcyc}
For cyclic extended semi-modules we have $\varphi(a+h)=\varphi(a)+1$ for all $a\in A$. Thus
\begin{equation*}
\mathcal{V}(A,\varphi)=\{(b_i,b)\mid b_i\in B, b\in A, b>b_i, \varphi(b)<\varphi(b_i)\}
\end{equation*}
where $B=A\setminus(h+A)$.
\end{remark}

\begin{prop}\label{propdim}
Let $(A,\varphi)$ be the cyclic extended semi-module associated to the normalized semi-module $A$ of type $\mu$. Then $\mid\mathcal{V}(A,\varphi)\mid = d(b,\mu)$.
\end{prop}
\begin{proof}
Recall that by $b_0$ we denote the minimal element of $A$ or $B$. Let $b_i$ be as in the definition of the type of $A$ and let $b_h=b_0$. First we show that 
\begin{eqnarray*}
\mathcal{V}(A,\varphi)&\rightarrow &\mathbb{Z}\\
(b_i,b)&\mapsto&b-b_i+b_h
\end{eqnarray*}
induces a bijection between $\mathcal{V}(A,\varphi)$ and $\{a\notin A\mid a>b_h\}$. Let $b\in A$ for some $b>b_i$. Then $b-b_i+b_{i+1}\notin A$ if and only if $(b_i,b)\in\mathcal{V}(A,\varphi)$. Let $b_{i_0}=\max\{b_i\in B\}$. We have $b\in A$ for all $b\geq b_{i_0}$. Thus for every $b>b_h$ with $b\notin A$, there is an element $(b_i,b-b_h+b_i)\in\mathcal{V}(A,\varphi)$ for some $h>i\geq i_0$. Hence $\{a\notin A\mid a>b_h\}$ is in the image of the map. To show that it is injective and that its image is contained in $\{a\notin A\mid a>b_h\}$, it is enough to show that $(b_i,b)\in\mathcal{V}(A,\varphi)$ implies that $b-b_i+b_j\notin A$ for all $j\in \{i+1,\dotsc,h\}$. Indeed, this ensures that $(b_j,b-b_i+b_j)\notin \mathcal{V}(A,\varphi)$ for all such $j$ and that $b-b_i+b_h\notin A$. We write $b=b_l+\alpha h$ for some $l$ and $\alpha$. Recall that $\varphi(b_i)=\mu_{i+1}$. As $(b_i,b_l+\alpha h)\in\mathcal{V}(A,\varphi)$, we have $\mu_{l+1}+\alpha<\mu_{i+1}$. Especially, $l<i$. This implies $\mu_{l+1}+\dotsm+ \mu_{l+\beta}+\alpha<\mu_{i+1}+\dotsm+\mu_{i+\beta}$ for all $\beta\leq h-i$. Using the recurrence for the $b_j$, one sees that this implies $b-b_i+b_{i+\beta}\notin A$ for all $\beta\leq h-i$.

It remains to count the elements of $\{a\notin A\mid a>b_0\}$. As $h+A\subseteq A$, we have
\begin{align*}
\mid\{a\notin A\mid a>b_0\}\mid&= \left(\sum_{i=0}^{h-1}b_i-b_0-i\right)\cdot\frac{1}{h}.\\
\intertext{From the construction of $A$ from its type we obtain}
&=\left(\sum_{i=0}^{h-1}\sum_{j=1}^{i}(m-\mu_jh)-i\right)\cdot\frac{1}{h}\\
&=\left(\sum_{i=0}^{h-1}\sum_{j=1}^{i}\frac{m}{h}-\mu_j\right)-\frac{h-1}{2}\\
&=d(b,\mu).
\end{align*}
\end{proof}

\begin{thm} \label{thm3}
Let $(A,\varphi)$ be an extended semi-module for $\mu$. Then
$\mid\mathcal{V}(A,\varphi)\mid\leq d(b,\mu)$. 
\end{thm}

\begin{proof}[Proof of Theorem \ref{thm3} for cyclic extended semi-modules]
We write $B=\{b_0,\dotsc,b_{h-1}\}$ as in the definition of the type $\mu'$ of $A$. As the extended semi-module is assumed to be cyclic, $\mu'$ is a permutation of $\mu$. Using Remark \ref{remcyc} we see
\begin{eqnarray*}
\mid \mathcal{V}(A,\varphi)\mid&=&\mid \{(b_i,a)\in B\times A\mid a>b_i, \varphi(a)<\varphi(b_i)\}\mid\\
&=& \sum_{\{(b_i,b_j)\in B\times B\mid b_j>b_i, \mu'_{j+1}<\mu'_{i+1}\}}\mu'_{i+1}-\mu'_{j+1}\\
&&+\mid\{(b_i,b_j+\alpha h)\mid b_j<b_i<b_j+\alpha h, \mu'_{i+1}>\mu'_{j+1}+\alpha\}\mid.
\end{eqnarray*}
We refer to these two summands as $S_1$ and $S_2$.

Let $(\tilde{b}_0,\tilde{\mu}_1),\dotsc,(\tilde{b}_{h-1},\tilde{\mu}_h)$ be the set of pairs $(b_0,\mu_1'),\dotsc,(b_{h-1},\mu_h')$, but ordered by the size of $b_i$. That is, $\tilde{b}_i<\tilde{b}_{i+1}$ for all $i$. Let
\begin{align*}
f:B&\rightarrow B\\
b_i&\mapsto b_{i+1}=b_i+m-\mu'_{i+1}h
\end{align*} where we identify $b_h$ with $b_0$.
This defines a permutation of $B$. From the ordering of the $\tilde{b}_i$ we obtain $\sum_{i=0}^{i_0}f(\tilde{b}_i)\geq\sum_{i=0}^{i_0}\tilde{b}_i$ for all $i_0$. As $f(\tilde{b}_i)=\tilde{b}_i+m-\tilde{\mu}_{i+1}h$, this is equivalent to $\sum_{i=1}^{i_0+1}\tilde{\mu}_i\leq (i_0+1)\frac{m}{h}$ for all $i_0$. We thus have $\nu\preceq\tilde{\mu}\preceq \mu$.

Recall the interpretation of $d(b,\mu)$ from Remark \ref{remgp}. We show that $S_1$ is equal to the number of lattice points above $\mu$ and on or below $\tilde{\mu}$. The second summand $S_2$ will be less or equal to the number of lattice points above $\tilde{\mu}$ and below $\nu$. Then the theorem follows for cyclic extended semi-modules.

We have $S_1=\sum_{i<j}\max\{\tilde{\mu}_{i+1}-\tilde{\mu}_{j+1},0\}$. Consider this sum for any permutation $\tilde{\mu}$ of $\mu$. If we interchange two entries $\tilde{\mu}_i$ and $\tilde{\mu}_{i+1}$ with $\tilde{\mu}_i>\tilde{\mu}_{i+1}$, the sum is lessened by the difference of these two values. There are also exactly $\tilde{\mu}_i-\tilde{\mu}_{i+1}$ lattice points on or below $\tilde{\mu}$ and above the polygon corresponding to the permuted vector. If $\tilde{\mu}=\mu$, both $S_1$ and the number of lattice points above $\mu$ and on or below $\tilde{\mu}$ are $0$. Thus by induction $S_1$ is equal to the claimed number of lattice points.

The last step is to estimate $S_2$. It is enough to construct a decreasing sequence (with respect to $\preceq$) of $\psi^i\in\mathbb{Q}^h$ for $i=0,\dotsc,h-1$ with $\psi^0=\tilde{\mu}$ and $\psi^{h-1}=\nu$ such that the number of lattice points above $\psi^i$ and on or below $\psi^{i+1}$ is greater or equal to the number of pairs $(\tilde{b}_{i+1},\tilde{b}_j+\alpha h)$ contributing to $S_2$. Note that the $\psi^i$ will no longer be lattice polygons. Let $f_i:B\rightarrow B$ be defined as follows: For $j>i$ let $f_i(\tilde{b}_j)=f(\tilde{b}_j)$. Let $\{f_i(\tilde{b}_j)\mid 0\leq j\leq i\}$ be the set of $f(\tilde{b}_j)$, but sorted increasingly. Let $\psi^i=(\psi^i_j)$ be such that $f_i(\tilde{b}_j)=\tilde{b}_j+m-\psi^i_{j+1} h$, i.e. $$\psi^i_{j+1}=\frac{\tilde{b}_j+m-f_i(\tilde{b}_j)}{h}=\frac{m}{h}-\frac{f_i(\tilde{b}_j)-\tilde{b}_j}{h}.$$ Similarly as for $\nu\preceq\tilde{\mu}$ one can show that 
\begin{equation*}
\nu\preceq\psi^{i+1}\preceq\psi^i\preceq\tilde{\mu}
\end{equation*}
for all $i$. As $f_0=f$ and $f_{h-1}=\id$, we have $\psi^0=\tilde{\mu}$ and $\psi^{h-1}=\nu$. It remains to count the lattice points between $\psi^i$ and $\psi^{i+1}$. To pass from $f_i$ to $f_{i+1}$ we have to interchange the value $f(\tilde{b}_{i+1})$ with all larger $f_i(\tilde{b}_j)$ with $j\leq i$. Thus to pass from the polygon associated to $\psi^i$ to the polygon of $\psi^{i+1}$ we have to change the value at $j$ by $(f_i(\tilde{b}_{j})-f(\tilde{b}_{i+1}))/h$, and that for all $j\leq i$ with $f_i(\tilde{b}_j)>f(\tilde{b}_{i+1})$. Thus there are at least 
$$\sum_{j\leq i, f_i(\tilde{b}_j)>f(\tilde{b}_{i+1})}\left\lfloor \frac{f_i(\tilde{b}_{j})-f(\tilde{b}_{i+1})}{h}\right\rfloor=\sum_{j\leq i, f(\tilde{b}_j)>f(\tilde{b}_{i+1})}\left\lfloor \frac{f(\tilde{b}_{j})-f(\tilde{b}_{i+1})}{h}\right\rfloor$$ lattice points above $\psi^i$ and on or below $\psi^{i+1}$. For fixed $i$ and $j<i+1$, the set of pairs 
$(\tilde{b}_{i+1},\tilde{b}_j+\alpha h)$ contributing to $S_2$ is in bijection with $\{\alpha\geq 1\mid f(\tilde{b}_j)-\alpha h>f(\tilde{b}_{i+1})\}$. The cardinality of this set is at most $\lfloor \frac{f(\tilde{b}_j)-f(\tilde{b}_{i+1})}{h}\rfloor$ which proves that $S_2$ is not greater than the number of lattice points between $\tilde{\mu}$ and $\nu$.
\end{proof}

\begin{ex}\label{ex4}
We give an example of a cyclic semi-module $(A,\varphi)$ where the type of $A$ is not dominant but where $\mid \mathcal{V}(A,\varphi)\mid =d(b,\mu)$. Let $m=4$, $h=5$, and $\mu=(0,0,1,1,2)$. Let $(A,\varphi)$ be the cyclic extended semi-module associated to the normalized semi-module of type $(0,0,1,2,1)$. Note that $A$ is the same semi-module as in Example \ref{ex2}. Then the dimension of the corresponding subscheme is 
$$\mid \mathcal{V}(A,\varphi)\mid =\mid \{(-1,2),(5,6),(5,7)\}\mid=d(b,\mu).$$
\end{ex}

\begin{proof}[Proof of Theorem \ref{thm3}]
Let $(A,\varphi)$ be an extended semi-module for $\mu$. Let $\varphi_i$ and $\mu^i$ be the sequences constructed in the proof of Lemma \ref{lemnoncyclic}. By induction on $i$ we show that $\mid\mathcal{V}(A,\varphi_i)\mid\leq d(b,\mu^i)$. For $i=0$, the extended semi-module $(A,\varphi_0)$ is cyclic, hence the assertion is already shown.

We use the notation of the proof of Lemma \ref{lemnoncyclic}. The description of the difference between $\mu^{i}$ and $\mu^{i-1}$ given there shows that
\begin{eqnarray*}
d(b,\mu^{i})-d(b,\mu^{i-1})&=&\sum_{l=1}^h\sum_{j=1}^l (\mu^{i-1}_{\dom,j}-\mu^{i}_{\dom,j})\\
\label{gldelta}&=&(\mid\{\mu^{i-1}_j\in (\varphi_{i-1}(x_i)-\alpha_i-n_i,\varphi_{i-1}(x_i)+1)\}\mid-1)\cdot\min\{\alpha_i,n_i+1\}.
\end{eqnarray*}
We denote this difference by $\Delta$. To show that $\mid\mathcal{V}(A,\varphi_{i})\mid-\mid\mathcal{V}(A,\varphi_{i-1})\mid\leq \Delta$ we use the decomposition into sequences $a^l_j$ of the extended semi-module $(A,\varphi_{i-1})$. Using the definition of $\mathcal{V}(A,\varphi)$ and the description of the difference between $\varphi_{i}$ and $\varphi_{i-1}$ from the proof of Lemma \ref{lemnoncyclic} one obtains
\begin{equation*}
\mid\mathcal{V}(a,\varphi_{i})\mid-\mid\mathcal{V}(a,\varphi_{i-1})\mid= S_1+S_2+S_3
\end{equation*}
where
\begin{eqnarray*}
S_1&=&\mid \{(x_i+h,b)\mid b\in A, b>x_i+h,\varphi_{i-1}(x_i)+1>\varphi_{i-1}(b)>\varphi_{i-1}(x_i)-\alpha_i\}\mid\\
S_2&=&\mid\{(b,x_i-\delta h)\mid b\in B\setminus\{x_i-n_ih\},b<x_i-\delta h, \delta\in \{0,\dotsc,n_i\},\\
&&\hspace{6cm}\varphi_{i-1}(x_i)-\delta-\alpha_i<\varphi_{i-1}(b)\leq \varphi_{i-1}(x_i)-\delta\}\mid\\
S_3&=&-\mid\{(x_i-n_ih,b)\mid b>x_i-n_ih,\varphi_{i-1}(x_i)-n_i>\varphi_{i-1}(b)\geq \varphi_{i-1}(x_i)-n_i-\alpha_i\}\mid.
\end{eqnarray*}
Here we used that $a\leq x_i$ implies that $\varphi_{i-1}(a+h)=\varphi_{i-1}(a)+1$. For each sequence $a^l_j$ of the extended semi-module $(A,\varphi_{i-1})$ we use $S_{1,l}$, $S_{2,l}$, and $S_{3,l}$ for the contributions of pairs with $b\in \{a^l_j\}$ to the three summands. Furthermore we write $S^l=S_{1,l}+S_{2,l}+S_{3,l}$. We show the following
assertions: If $\varphi_{i-1}(a^l_0)\notin (\varphi_{i-1}(x_i)-\alpha_i-n_i,\varphi_{i-1}(x_i)+1)$ or if $a^l_0=x_i-n_ih$, then $S^l=0$. Otherwise, $S^l\leq \min\{\alpha_i,n_i+1\}$. Then the theorem follows from property (4c) of extended semi-modules.

To determine the $S^l$, we consider the following cases:

{\bf Case 1:} $\varphi_{i-1}(a_0^l)\geq \varphi_{i-1}(x_i)+1$. In this case it is easy to see that $S_{1,l}=S_{2,l}=S_{3,l}=0$.

{\bf Case 2:} $a_0^l>x_i$. This implies that $S_{2,l}=0$. If $\varphi_{i-1}(a_0^l)\leq \varphi_{i-1}(x_i)-n_i-\alpha_i$, then $S_{1,l}+S_{3,l}=\alpha_i-\alpha_i=0$. 
Let now $\varphi_{i-1}(a_0^l)\in (\varphi_{i-1}(x_i)-\alpha_i-n_i,\varphi_{i-1}(
x_i)+1)$. Then
\begin{equation*}
S_{1,l}+S_{3,l}\leq \mid\{ a_j^l\mid \varphi_{i-1}(x_i)+1>\varphi_{i-1}(a_j^l)\geq \max\{\varphi_{i-1}(x_i)-\alpha_i+1,\varphi_{i-1}(x_i)-n_i\}\}\mid.
\end{equation*}
As $\varphi_{i-1}(a_{j+1}^l)=\varphi_{i-1}(a_j^l)+1$ for all $j$, the right hand side is less or equal to $\min\{\alpha_i,n_i+1\}$.

{\bf Case 3:} $a_0^l=x_i-n_ih$. This sequence starts with $x_i-n_ih,\dotsc, x_i,x_i+h$. (Recall that the sequences $\{a_j^l\}$ for $\varphi_{i-1}$ are of this easy form with stepwidth $h$ as long as $a_j^l\leq x_i<x_{i-1}$.) Note that within one sequence $a_j^l>a_{j'}^l$ implies $\varphi_{i-1}(a_j^l)>\varphi_{i-1}(a_{j'}^l)$. Hence this special sequence does not make any contribution, as in $S^l$ we only consider pairs where both elements are in the sequence starting with $x_i-n_ih$.

{\bf Case 4:} $a_0^l<x_i$, but not congruent to $x_i$ modulo $h$. Again $a_{j+1}^l=a_j^l+h$ if $a_j^l\leq x_i$. We first assume that $\varphi_{i-1}(a_0^l)\leq \varphi_{i-1}(x_i)-n_i-\alpha_i$. Then $S_{2,l}=0$. Assume that $b=a_j^l$ contributes to $S_{1,l}$. Then $j\geq n_i+1$ and $a_j^l>x_i+h$. If $a_0^l<x_i-n_ih$, then $[x_i-n_ih, x_i+h]$ contains $n_i+1$ elements of the sequence. Thus in all cases $a^l_{j-n_i-1}>x_i-n_ih$. This element then leads to a contribution to $S_{3,l}$, as $\varphi_{i-1}(a^l_{j-n_i-1})=\varphi_{i-1}(a^l_j)-n_i-1$. In the other direction, if $a^l_j$ contributes to $S_{3,l}$, then $a^l_{j+n_i+1}$ contributes to $S_{1,l}$. Thus $S^l=0$. We now assume that $\varphi_{i-1}(a_0^l)\in (\varphi_{i-1}(x_i)-\alpha_i-n_i,\varphi_{i-1}(x_i)+1)$. Let $n$ be maximal with $a_n^l=a_0^l+nh<x_i$. Then we have
\begin{eqnarray*}
S_{1,l}&=&\mid\{a_j^l\mid j>n+1,\varphi_{i-1}(x_i)\geq\varphi_{i-1}(a_0^l)+j>\varphi_{i-1}(x_i)-\alpha_i\}\mid\\
S_{2,l}&=&\mid\{a_j^l\mid 0\leq j\leq \min\{n,n_i\},\varphi_{i-1}(x_i)\geq\varphi_{i-1}(a_0^l)+j>\varphi_{i-1}(x_i)-\alpha_i\}\mid\\
S_{3,l}&=&-\mid\{a_j^l\mid j\geq\max\{n-n_i+1,0\},\varphi_{i-1}(x_i)-n_i>\varphi_{i-1}(a_0^l)+j\}\mid \\
&=&-\mid\{a_j^l\mid j>\max\{n+1,n_i\},\varphi_{i-1}(x_i)\geq\varphi_{i-1}(a_0^l)+j\}\mid.
\end{eqnarray*}
Thus $$S^l\leq S_{1,l}+S_{2,l}\leq \{j\mid \varphi_{i-1}(x_i)\geq\varphi_{i-1}(a_0^l)+j>\varphi_{i-1}(x_i)-\alpha_i\}=\alpha_i.$$
If $n+1\geq n_i$, then $S_{1,l}+S_{3,l}\leq 0$. Thus $S^l\leq S_{2,l}\leq n_i+1$. If $n_i>n+1$ then $S_{1,l}+S_{3,l}\leq n_i-n-1$ and $S_{2,l}\leq n+1$. Hence in both cases $S^l\leq\min\{\alpha_i,n_i+1\}$.
\end{proof}
\begin{ex}\label{ex3}
Example \ref{ex2} describes a non-cyclic extended semi-module $(A,\varphi)$ for $\mu=(0,0,0,2,2)$  such that   
$$\mid \mathcal{V}(A,\varphi)\mid =\mid \{(5,6),(5,7),(4,6),(4,7)\}\mid=d(b,\mu).$$
\end{ex}

\begin{proof}[Proof of Theorem \ref{thm1}]
Lemma \ref{lemdim} and Theorem \ref{thm2} imply that $\dim X_{\mu}(b)^0=\max\mid \mathcal{V}(A,\varphi)\mid$. In Proposition \ref{propdim} we give a pair with $\mid\mathcal{V}(A,\varphi)\mid=d(b,\mu)$. Theorem \ref{thm3} shows that the maximum is at most $d(b,\mu)$. Together we obtain $\dim X_{\mu}(b)=d(b,\mu)$.
\end{proof}

\section{Irreducible components}
\begin{kor}
Let $G=GL_h$, let $b$ be superbasic and $\nu\preceq\mu$. Then the action of $J(F)$ on the set of irreducible components of $X_{\mu}(b)$ has only finitely many orbits.
\end{kor}
\begin{proof}
It is enough to consider the intersection of the orbits with the set of irreducible components of $X_{\mu}(b)^0$. Theorem \ref{thm2} implies that each $\mathcal{S}_{A,\varphi}$ is irreducible. Thus the Corollary follows from Lemma \ref{lemmufin}. 
\end{proof}
\begin{ex}
We give two examples to show that even for superbasic $b$, the irreducible components of $X_{\mu}(b)$ are in general not permuted transitively by $J(F)$. The description of $J(F)$ in Section \ref{section2} implies that $A(gM)=A(M)$ and $\varphi(gM)=\varphi(M)$ for each $g\in J(F)$ with $v_t(\det(g))=0$. First 
we consider the example $m=4$, $h=5$, and $\mu=(0,0,1,1,2)$. It is enough to find two extended semi-modules for $\mu$ leading to subschemes of dimension $d(b,\mu)=3$. Indeed, the subschemes corresponding to different extended semi-modules are disjoint and lead to irreducible components in different $J(F)$-orbits. One such extended semi-module is the cyclic extended semi-module considered in Proposition \ref{propdim}. A second extended semi-module $(A,\varphi)$ is given in Example \ref{ex4}. Here, $A$ is of type $(0,0,1,2,1)$, hence different from the semi-module considered before.

For the second example let $m=4$, $h=5$, and $\mu=(0,0,0,2,2)$. Here the two extended semi-modules for $\mu$ leading to subschemes of dimension $d(b,\mu)=4$ are the ones considered in Proposition \ref{propdim} and Examples \ref{ex2} and \ref{ex3}. The corresponding semi-modules are different as they are of type $(0,0,0,2,2)$ and $(0,0,1,2,1)$.
\end{ex}

\end{document}